\documentclass[openright,a4paper,american,11pt]{article}

\usepackage[latin1]{inputenc}

\usepackage{amsmath}

\usepackage{babel}

\usepackage{amsfonts}
\input xy
\xyoption{all}
\usepackage[dvips]{graphicx}
\usepackage{boxedminipage}
\usepackage{color}

\newtheorem{thm}{Theorem}[section]
\newtheorem{defi}[thm]{Definition}
\newtheorem{prop}[thm]{Proposition}
\newtheorem{lemma}[thm]{Lemma}
\newtheorem{cor}[thm]{Corollary}
\newcommand{\findemo}{\hfill
                     $\Box$ \vspace{1.5 ex}}

\begin{document}
\title{REAL PLANE ALGEBRAIC CURVES WITH ASYMPTOTICALLY MAXIMAL NUMBER OF EVEN OVALS}
\author{Erwan Brugall\'e}
\date{}
\maketitle
\begin{abstract}
It is known for a long time that a nonsingular real algebraic curve of degree $2k$ in the projective plane cannot have more than $\frac{7k^2}{2}-\frac{9k}{4}+\frac{3}{2}$ even ovals. We show here that this upper bound is asymptotically sharp, that is to say we construct a family of curves of degree $2k$ such that $\frac{p}{k^2}\to_{k\to\infty} \frac{7}{4}$, where $p$ is the number of even ovals of the curves. We also show that the same kind of result is valid dealing with odd ovals.
\end{abstract}

\section{Introduction}
The set $\mathbb RA$ of the real points of a nonsingular real
algebraic curve in $\mathbb RP^2$ is a disjoint union of circles. Each
of these circles either disconnects $\mathbb RP^2$ or not. In the
former case, such a circle is called \textit{an oval} of $\mathbb
RA$. The part of $\mathbb RP^2$ cut along this oval which is
homeomorphic to a disk is called \textit{the interior} of the oval. An
oval of a curve which contains no other ovals of the curve is called
\textit{an empty oval}. All the connected components of a curve of
even degree are ovals. For such a curve, the ovals which are contained
in an even (resp. odd) number of ovals are called \textit{even}
(resp. \textit{odd}) ovals. Given a real plane algebraic curve of
degree $2k$, its number of even (resp. odd) ovals will be denoted by
$p$ (resp. $n$). This separation of ovals in two groups is important
for many reasons. One of them is the fact that
curves with many even ovals can be used to construct real
algebraic surfaces with big Betti numbers (see section
\ref{surfaces}).

\noindent What are the maximal possible values for $p$ and $n$ with
respect to $k$? The first step in the study of this problem is due to
V. Ragsdale. In 1906, she conjectured in \cite{Rag} that
$p\le\frac{3k(k-1)}{2}+1$ and $n\le\frac{3k(k-1)}{2}$. About 30 years
later, I. G. Petrovsky proved in \cite{Pet} that
$p-n\le\frac{3k(k-1)}{2}+1$ and $n-p\le\frac{3k(k-1)}{2}$ (these
inequalities were also conjectured by Ragsdale), and formulated a
conjecture  similar to Ragsdale's one (it seems clear that Petrovsky
was not familiar with  Ragsdale's work). Combining the first Petrovsky
inequality with the  Harnack Theorem (which gives an upper bound for
the number of connected  components of a real plane algebraic curve
with respect to its degree)   one can obtain the following upper
bounds for $p$ and $n$ : 
$$p\le \frac{7k^2}{4}-\frac{9k}{4}+\frac{3}{2}\textrm{   and    } n\le \frac{7k^2}{4}-\frac{9k}{4}+1.$$
The first counterexamples to Ragsdale's conjecture for $n$ (but not to
Petrovsky's one) were constructed by O. Ya. Viro in the late 70's (see
\cite{V4}). In 1993, I. Itenberg gave in \cite{I3} counterexamples to
Ragsdale's and Petrovsky's conjectures. He has constructed  for every
positive integer $k$ curves of degree $2k$ with $\frac{13k^2}{8}+O(k)$
even ovals and curves of degree $2k$ with $\frac{13k^2}{8}+O(k)$ odd
ovals. These lower bounds were successively improved by B. Haas (see
\cite{Haa}) and Itenberg (see \cite{I1}). The best lower bound known
before the present paper was $\frac{81k^2}{48}+O(k)$ for both $p$ and
$n$. We point out the fact that that no counterexamples of Ragsdale's conjectures is known among curves with the maximal number of connected components.

\noindent All these constructions are based on a particular case of the so called Viro method (see \cite{V1}, \cite{V2}), the combinatorial patchworking. One can note that dealing with non convex triangulations (and so with pseudo-holomorphic curves, see \cite{IS}), F. Santos (\cite{San}) constructed curves with $\frac{17k^2}{10}+O(k^{\frac{3}{2}})$ even ovals.

\noindent It seemed to us that the $T$-construction is more or less
``rigid'' and that the general Viro method gives  one more flexibility
and possibilities to construct real algebraic curves. Then, we resumed
the work of Itenberg and Santos in this way, trying to increase the
density of even ovals. It turned out that gluing curves whose Newton
polygon is not anymore a triangle but an hexagon, it was possible to
prove that the upper bounds given by the
Harnack theorem
and the Petrovsky
inequalities
are asymptotically
sharp. 

This is the main result of this article.

\begin{figure}[h]
\centering
\begin{tabular}{c}
\includegraphics[height=15cm, angle=0]{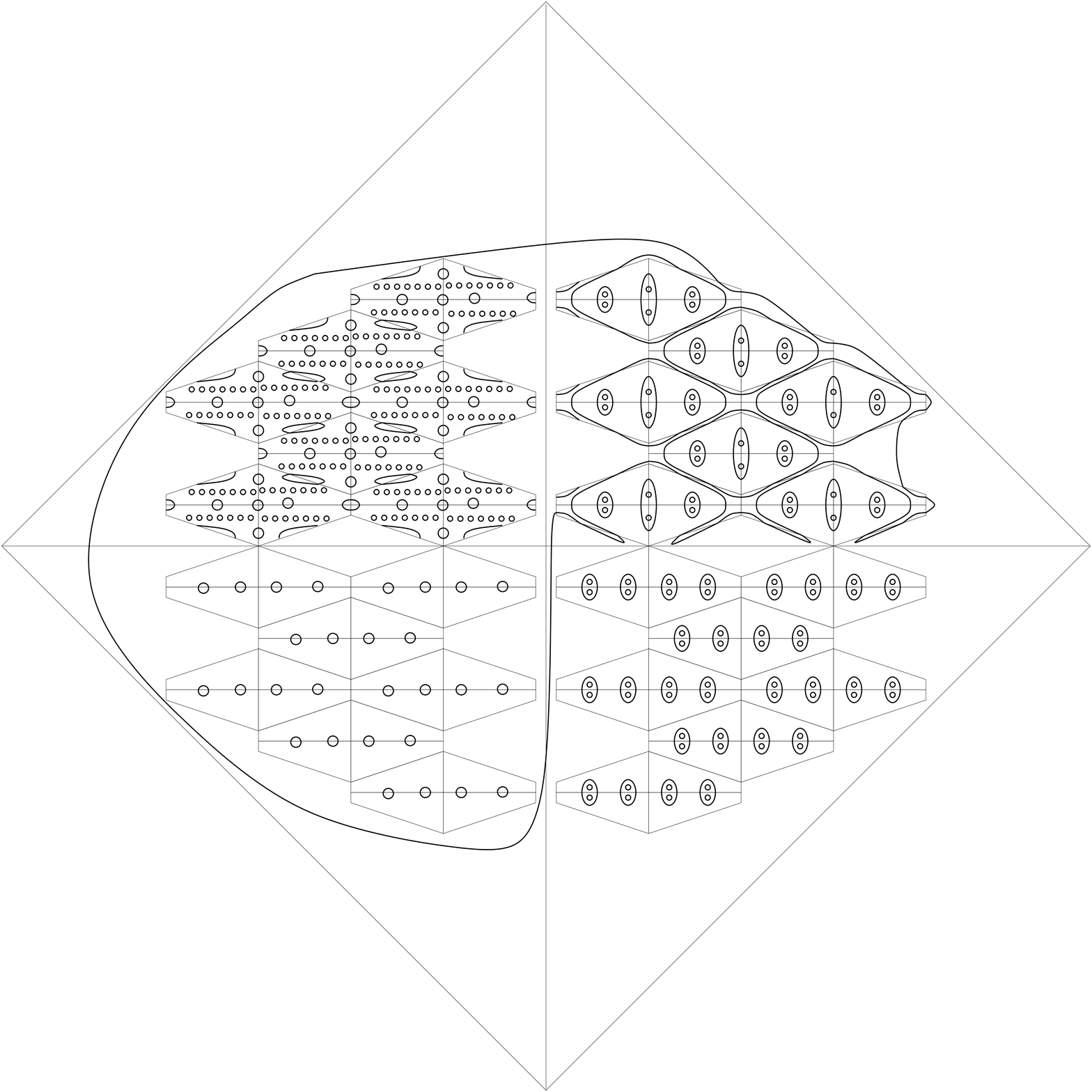}
\\ 
\end{tabular}
\caption{}
\label{Viro big n}
\end{figure}

\begin{thm}\label{max even}
\noindent There exists a family $(C_{2k})_{k\ge 0}$ of nonsingular
real algebraic curves of degree $2k$ in $\mathbb RP^2$ such that
$$\lim_{k\to\infty} \frac{p}{k^2}=\frac{7}{4}. $$

\noindent There exists a family $(C_{2k})_{k\ge 0}$ of nonsingular
real algebraic curves of degree $2k$ in $\mathbb RP^2$ such that
$$\lim_{k\to\infty} \frac{p}{k^2}=\frac{7}{4}. $$
\end{thm}
\textit{Proof. }The assertion relative to $p$ is a direct consequence of corollary \ref{many even}. The assertion relative to $n$ can be proved, as in \cite{I3} and \cite{I1}, by a small modification of the construction given in section \ref{construct} (see Figure \ref{Viro big n}).\findemo

This article is organized as follows : in section \ref{defi}, we
recall some facts about rational geometrically ruled surfaces. In
section \ref{construct}, we prove the first part of Theorem \ref{max
  even}. The constructions in this section are based on the classical
Viro method. We assume in this section the existence of some special curves in rational geometrically ruled surfaces. The construction of the latter curves are based on the less classical \textit{real rational graphs} theoretical method. Section \ref{graph} is devoted to the definition and properties of such graphs, and in section \ref{Constr tangency}, we construct the special curves used in section \ref{construct}. In section \ref{surfaces}, we give some applications of Theorem \ref{max even} to real algebraic surfaces.

\vspace{2ex}
\noindent \textbf{Acknowledgment. }I am grateful to Ilia Itenberg for useful discussions and advisements and to Olivier Le Gal for his art in hexagons counting.

\section{Rational geometrically ruled surfaces}\label{defi}
\textit{The $n^{th}$ rational geometrically ruled surface}, denoted by $\Sigma_n$, is the surface obtained by taking four copies of $\mathbb C^2$ with coordinates $(x,y)$, $(x_2,y_2)$, $(x_3,y_3)$
and $(x_4,y_4)$, and by gluing them along $(\mathbb C^*)^2$ with the identifications $(x_2,y_2)=(1/x,y/x^n)$, $(x_3,y_3)=(x,1/y)$
and $(x_4,y_4)=(1/x,x^n/y)$. Let us denote by $E$  (resp. $B$ and $F$) the algebraic curve in $\Sigma_n$ defined by the equation $\{y_3=0\}$ (resp. $\{y=0\}$ and $\{x=0\}$). The coordinate system $(x,y)$
is called \textit{standard}.
The projection $\pi$ : $(x,y)\mapsto x$ on $\Sigma_n$ defines a $\mathbb CP^1$-bundle over $\mathbb CP^1$.
We have $B\circ B=n$, $F\circ F=0$ and $B\circ F=1$. The surface $\Sigma_n$ has a natural real structure induced by the complex conjugation in $\mathbb C^2$, and
the real part of $\Sigma_n$ is a torus if $n$ is even and a Klein bottle if $n$ is odd. The restriction of $\pi$ on $\mathbb R\Sigma_n$ defines a pencil of lines denoted by $\mathcal L$.

 The group $H_2(\Sigma_n,\mathbb Z)$ is isomorphic to $\mathbb Z\times\mathbb Z$ and is generated by the classes of $B$ and $F$. Moreover, one has $E=B-nF$. An algebraic curve on $\Sigma_n$ is
said to be of \textit{bidegree} $(k,l)$ if it realizes the homology class $kB+lF$ in $H_2(\Sigma_n,\mathbb Z)$. A curve of bidegree $(3,0)$ is called \textit{a trigonal curve} on $\Sigma_n$.

In the rational geometrically ruled surfaces, we study real curves up to \textit{isotopy with respect to $\mathcal L$}. Two curves are said
to be isotopic with respect to the fibration $\mathcal L$ if there exists an isotopy of $\Sigma_n$ which brings the first curve to the second one, and
which maps each line of $\mathcal L$ in another line of $\mathcal L$.

\noindent In this paper, curves in a rational geometrically ruled surface are depicted up to isotopy with respect to $\mathcal L$.

\section{Construction of real algebraic curves with many even ovals}\label{construct}
We will use the Viro method to construct real plane algebraic curves.  The unfamiliar readers can refer to \cite{V1} and \cite{V2}.

\begin{figure}[h]
\centering
\begin{tabular}{cc}
\includegraphics[height=3cm, angle=0]{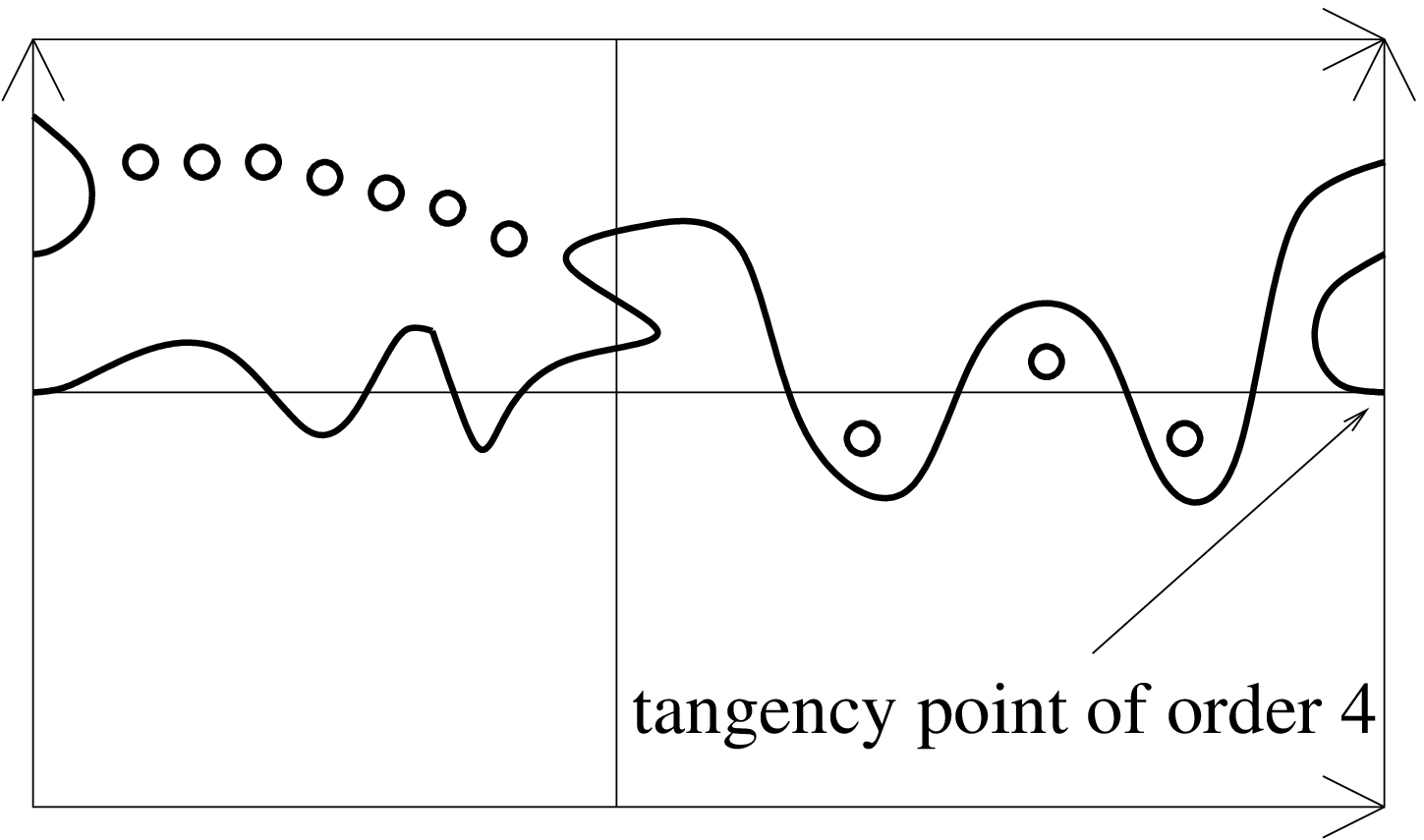}&
\includegraphics[height=3cm, angle=0]{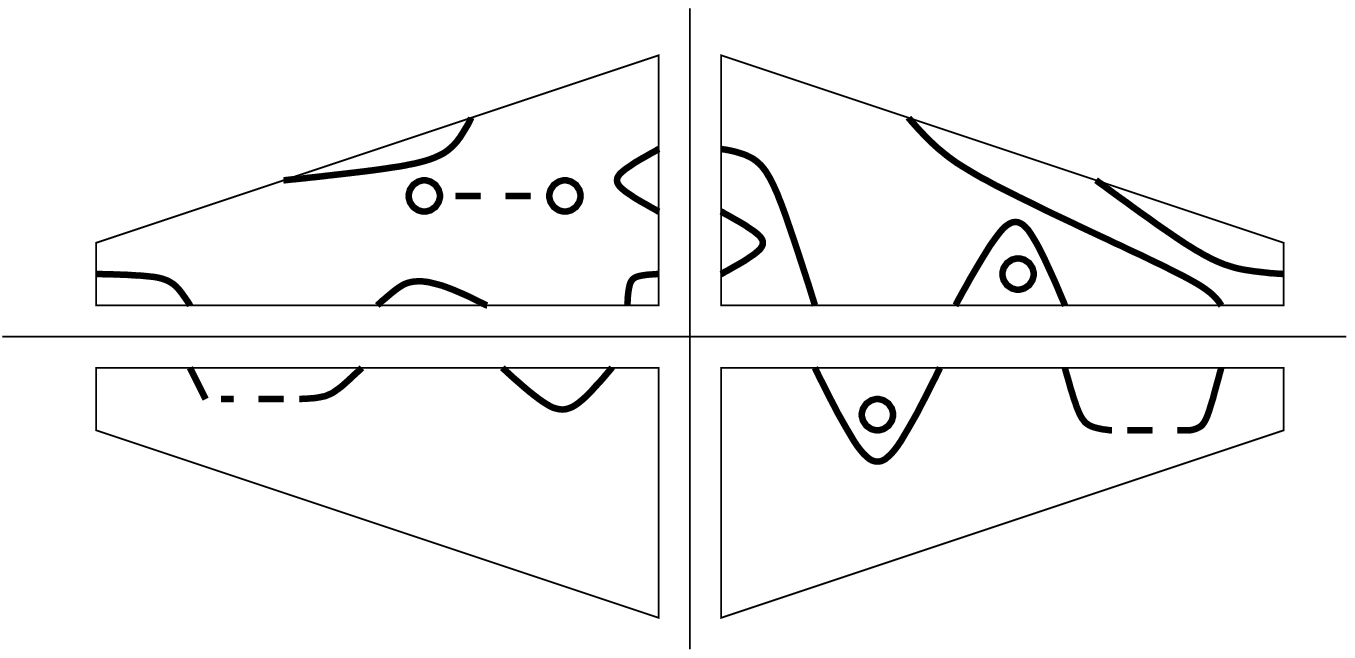}
\\ a)&b)
\end{tabular}
\caption{}
\label{ex n=4}
\end{figure}

\noindent In this section, we will use the following proposition which will be proved in section \ref{Constr tangency}.
\begin{prop}\label{C_n}
\begin{figure}[h]
\centering
\begin{tabular}{cc}
\includegraphics[height=3cm, angle=0]{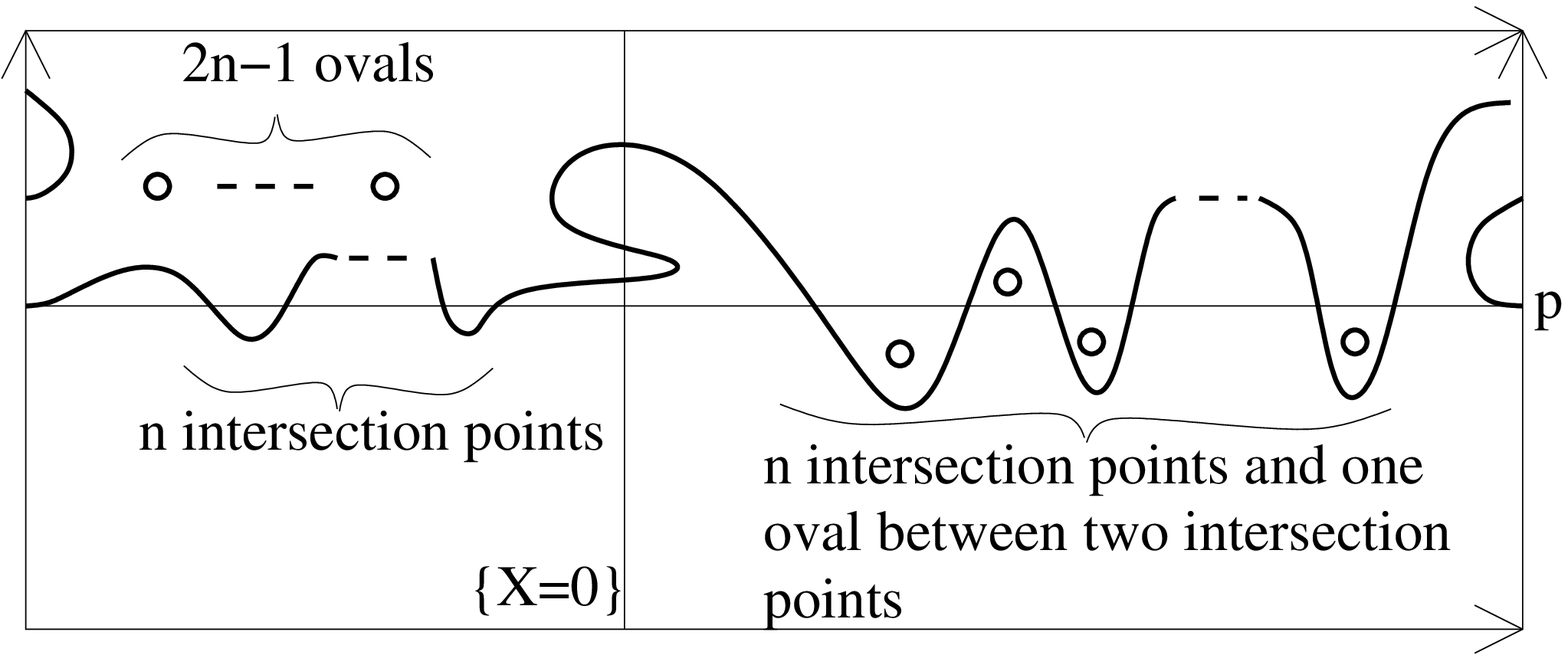}&
\includegraphics[height=3cm, angle=0]{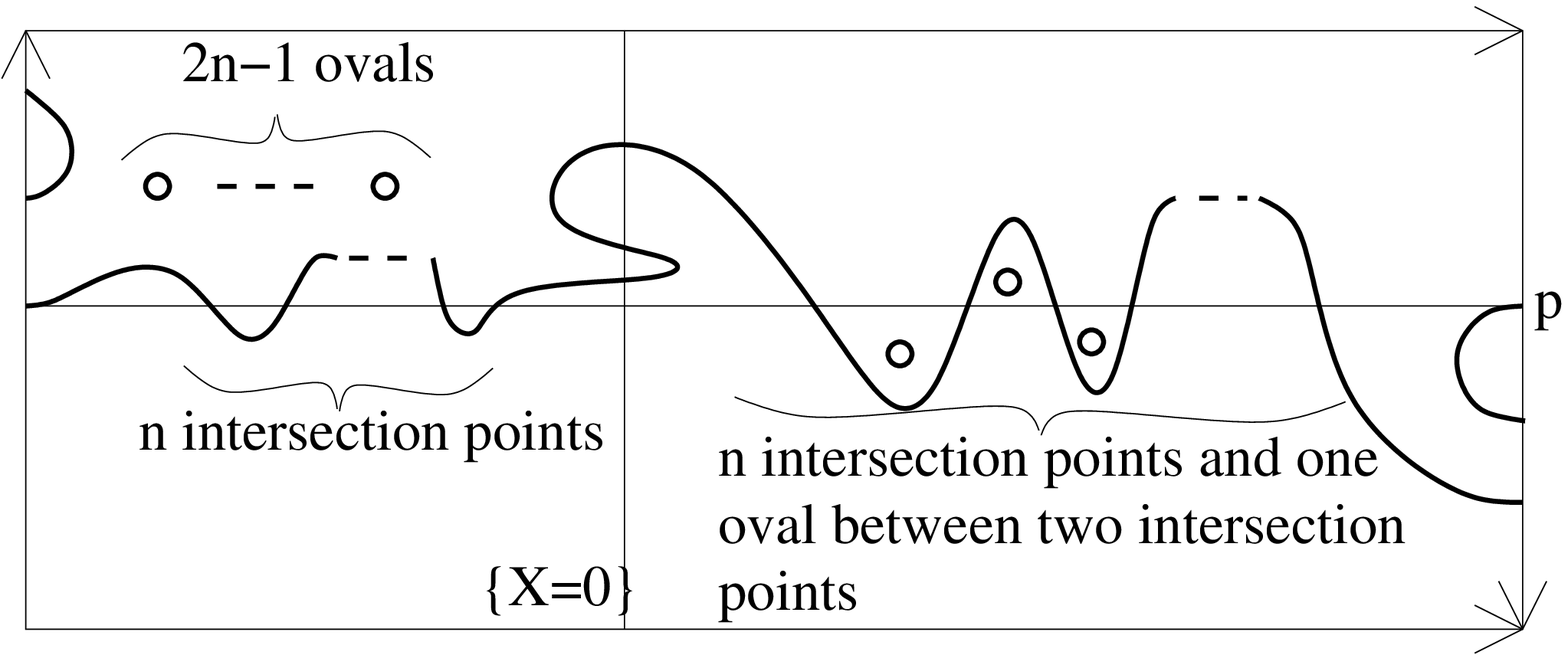}
\\ a)&b)
\end{tabular}
\caption{}
\label{Tangent curves}
\end{figure}
For any $n\ge 1$, there exists a maximal real algebraic trigonal curve $C_n$ in $\Sigma_n$ realizing the $\mathcal L-scheme$ and whose position with respect to the axis $\{Y=0\}$ and $\{X=0\}$ is depicted in Figure \ref{Tangent curves}a) if $n$ is even and \ref{Tangent curves}b) if $n$ is odd, where $p$ is a tangency point of order $n$ of $C_n$ and the axis $\{Y=0\}$.
\end{prop}
For example, the curve for $C_4$ is depicted on figure \ref{ex n=4}.

\vspace{2ex}
\noindent Let us fix an even integer $n$.

\begin{figure}[h]
\centering
\begin{tabular}{cc}
\includegraphics[height=5cm, angle=0]{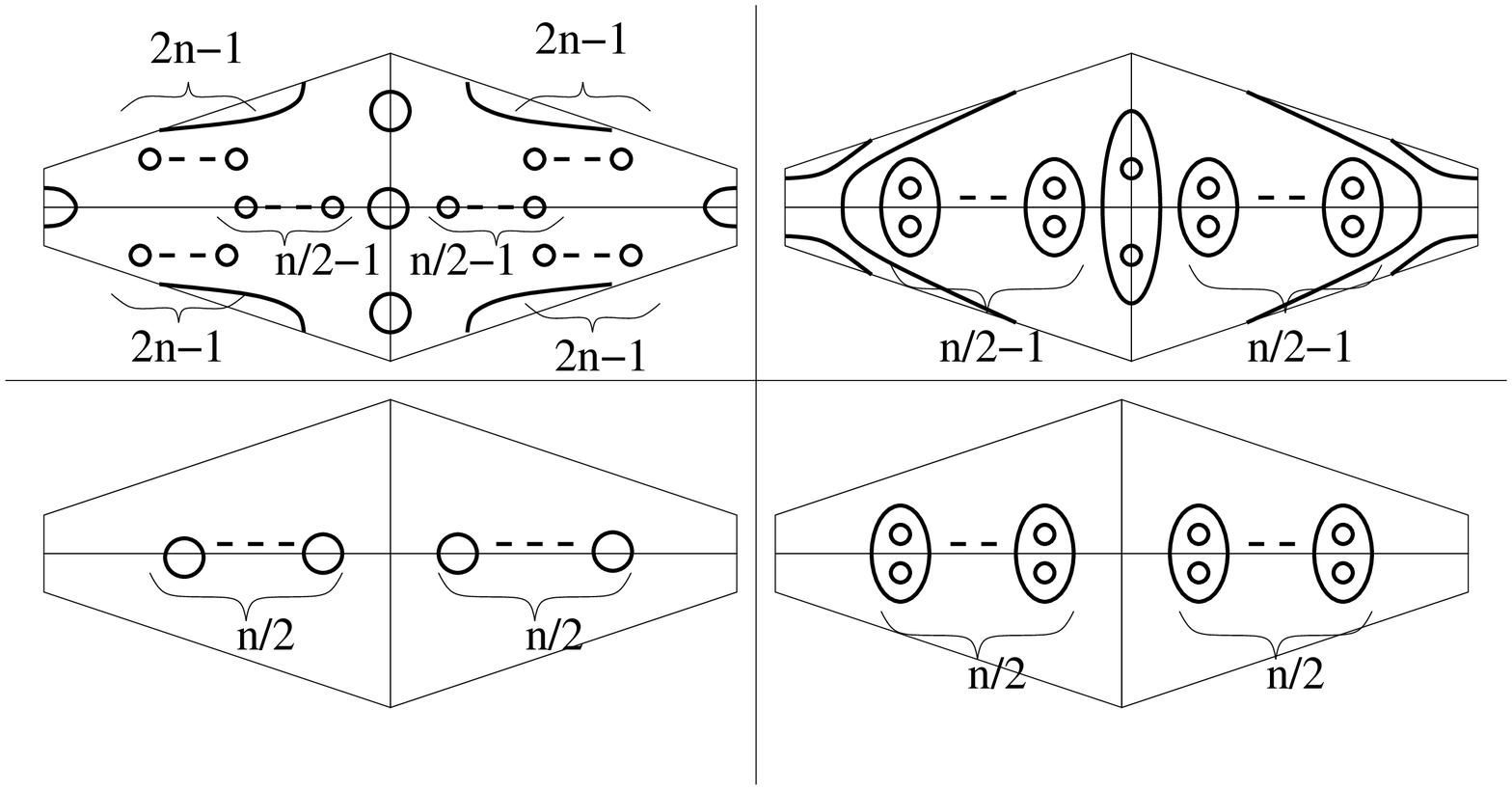}&
\includegraphics[height=5cm, angle=0]{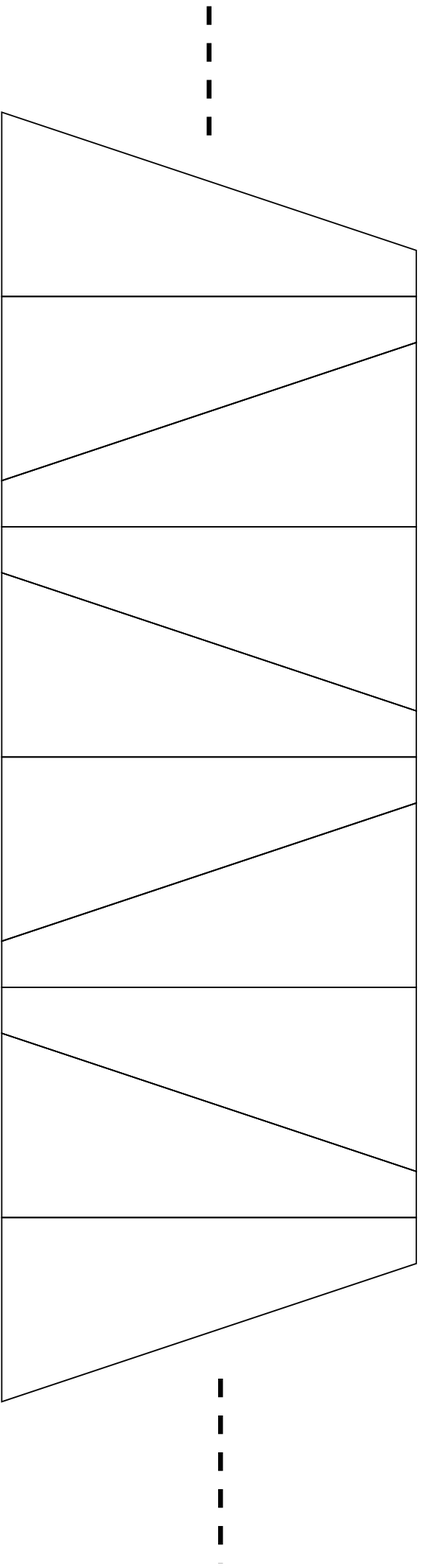}
\\ a)&b)
\end{tabular}
\caption{}
\label{chart C_n}
\end{figure}

\noindent The Newton polygon of the curve $C_n$ is the quadrangle with vertices $(0,0)$, $(2n,0)$, $(2n,1)$ and $(0,3)$ and the chart of $C_n$ is depicted on Figure \ref{ex n=4}b) (we disjointed the 4 symmetric copies of the Newton polygon of $C_n$ for convenience). Moreover, performing the transformation $\widetilde{Y}=\lambda Y$ if necessary, we can assume that the truncation of $C_n$ on  $[(0,3);(2n,1)]$ is $\alpha Y^3+\beta Y^2X^n+\alpha YX^{2n}$ with $\alpha$ and $\beta$ two real numbers.
Let us denote by $H_n$ the hexagon obtained by gluing the charts of 4 birational transforms of $C_n$ as depicted on Figure \ref{chart C_n}a).

\begin{figure}[h]
\centering
\begin{tabular}{c}
\includegraphics[height=15cm, angle=0]{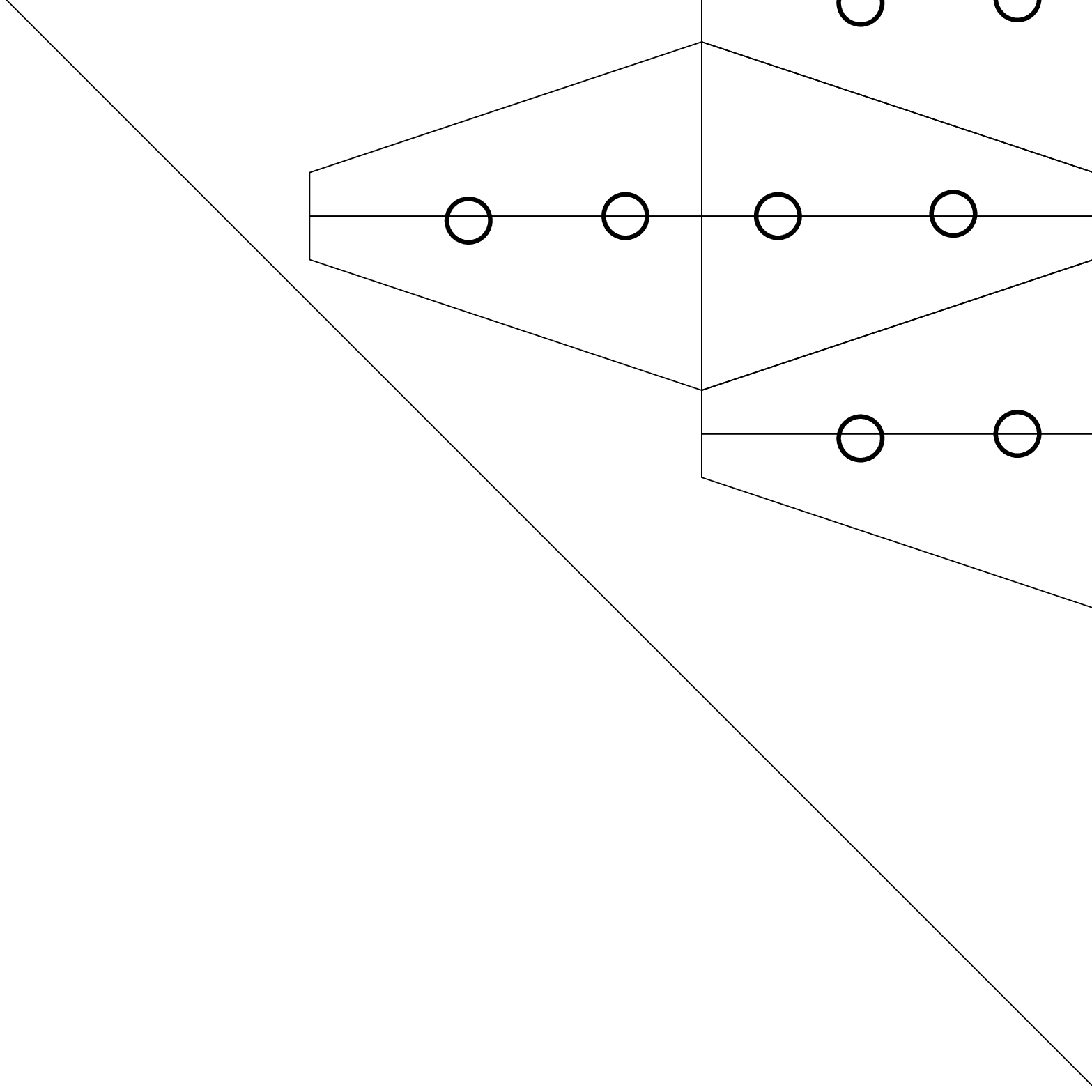}
\\ 
\end{tabular}
\caption{}
\label{begin Viro}
\end{figure}

\noindent Let us fix an integer $k$ and denote by $T_{2k}$ the triangle with vertices $(0,0)$, ($2k,0)$ and $(0,2k)$. We start a subdivision of $T_{2k}$ in the following way : for each integers $l$ and $h$, we put the hexagon $H_n$ centered in the point $(1+2n+4l,3+8h)$ or $(1+4n+4l,7+8h)$ if this hexagon is contained in $T_{2k}$. In this way, we obtain the beginning of a patchwork of a real plane curve of degree $2k$ as depicted on Figure \ref{begin Viro} (here were chose n=4 for convenience).
\begin{lemma}\label{convexity}
The union of all the hexagons can be a part of a convex subdivision of $T_{2k}$.
\end{lemma}
\textit{Proof. }The union of the hexagons can be decomposed in vertical strips as depicted on Figure \ref{chart C_n}b). Given any convex function on the left edge of the strip, one can extend it to a convex function on the whole strip which induces this subdivision.\findemo

\begin{figure}[h]
\centering
\begin{tabular}{c}
\includegraphics[height=15cm, angle=0]{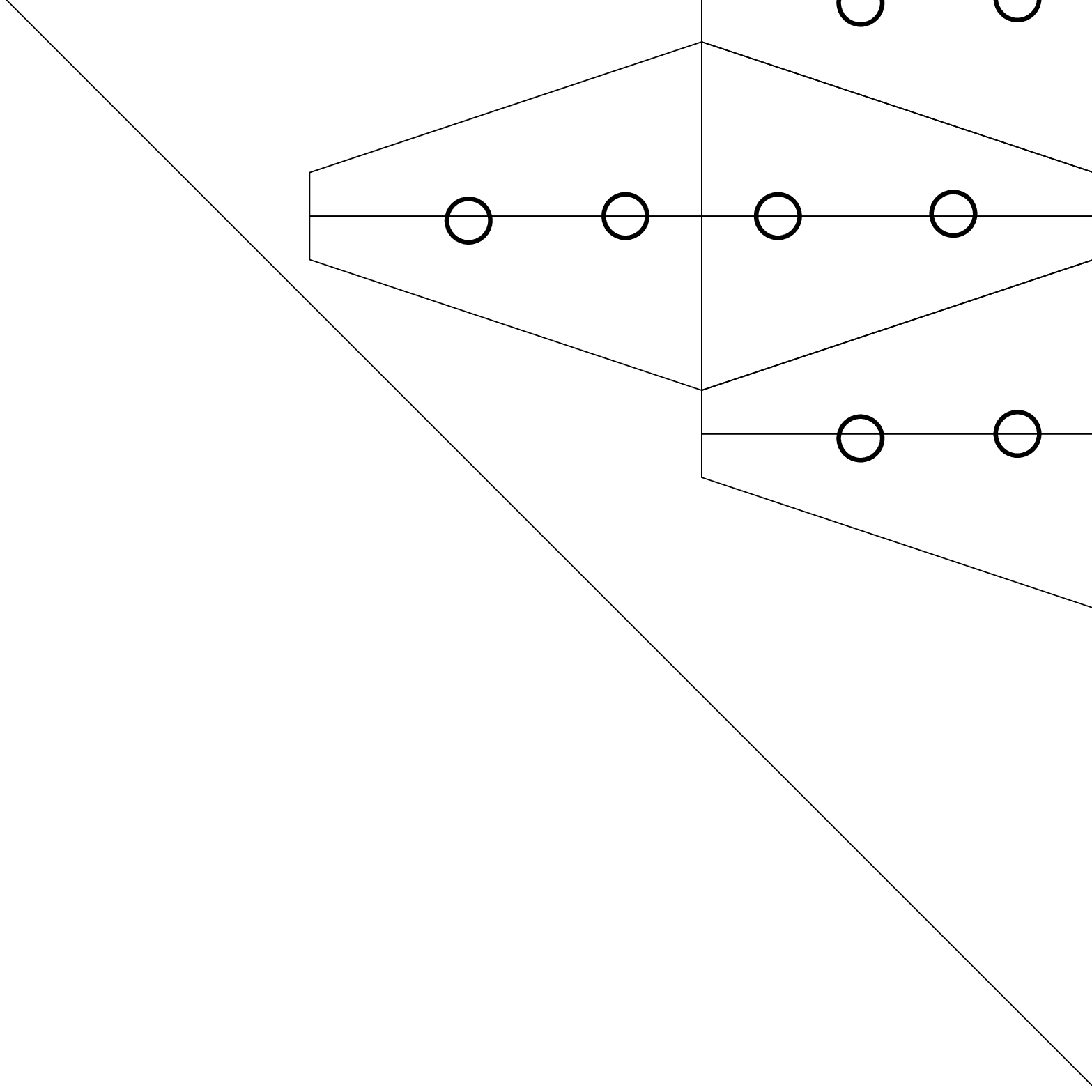}
\\ 
\end{tabular}
\caption{}
\label{finish Viro}
\end{figure}

\noindent Suppose we are given an extension of our beginning of patchwork to the whole $T_{2k}$, satisfying the hypothesis of the Viro Theorem (see \cite{V1}). Then, by Viro's Theorem, we obtain a real algebraic curve of degree $2k$ in $\mathbb RP^2$, which we will denote by $C_k^n$.
Next, let us choose such an extension such that (see Figure \ref{finish Viro}) :
\begin{itemize}
\item each oval of $C_k^n$ lying in the half plane $\{x<0\}$ and coming from an hexagon is even and not contained in another oval of the curve,
\item each oval of $C_k^n$ lying in the quadrant $\{x>0\}$ and coming from an empty oval of an hexagon is even and contained in two other ovals of the curve,
\end{itemize}
where the coordinates system is the one given by the chart of $C_k^n$.

\noindent It is clear  that such an extension exists. One can
construct curves we need to complete our patchwork   by
the classical small perturbation method (see, for example, \cite{V4}). The convexity condition can be ensured, for example, by keeping on decomposing $T_{2k}$ in strips.
\begin{lemma}\label{even H_n}
Each hexagon contributes of 
at least $14n-5$ even ovals to the curve $C_k^n$.
\end{lemma}
\textit{Proof. }Straightforward.\findemo

\begin{lemma}\label{lower bound}
The curve $C_k^n$ has at least $\frac{7k^2}{4}-\frac{5k^2}{8n}-21nk+\frac{15k}{2}+63n^2-\frac{45n}{2}$ even ovals.
\end{lemma}
\textit{Proof. }According to Lemma \ref{even H_n}, each hexagon $H_n$ in the patchwork of the curve $C_k^n$ gives at least $14n-5$ even ovals. Then, if the patchwork contains $N$ hexagons $H_n$, the curve  $C_k^n$ will have at least $N(14n-5)$ even ovals.
\noindent The triangle $T'_{2k}$ with vertices $(6n,6n)$, $(2k-12n,6n)$ and $(6n,2k-12n)$ is contained in the union  of the hexagons, so
$$N\ge \frac{Area(T'_{2k})}{Area(H_n)}=\frac{(k-6n)^2}{8n}.$$
Hence the number of even ovals of $C_k^n$ is at least $\frac{(k-6n)^2}{8n}(14n-5)$. Developing this quantity, we obtain the lower bound stated in the lemma.\findemo

\noindent The same construction can be done with an odd integer $n$. The curve obtained is also denoted by $C_k^n$ and the lower bound of lemma \ref{lower bound} for its number of even oval is still valid.

\noindent Now we are able to prove the main theorem of this paper. We denote the integer part of a real $r$ by $[r]$.
\begin{cor}\label{many even}
The curve $C_k^{[\sqrt{k}]}$ has $\frac{7k^2}{4}+O(k^{\frac{3}{2}})$
even ovals.\findemo 
\end{cor}

\section{Real rational graphs on $\mathbb CP^1$}\label{graph}
This section deals with the following problem : given a real
arrangement of roots of three real polynomials (called \textit{a root
scheme} below), does there exist two real polynomials $P$ and $Q$ such
that the real roots of $P,Q$ and $P+Q$ realize the given arrangement?

\noindent This question can be reformulated in terms of existence of a
certain graph on $\mathbb CP^1$ (called \textit{a real rational graph}
below). 

\noindent We start with the following fact :
to any
rational
map
$f:\mathbb CP^1\to\mathbb CP^1$, one can associate
a natural graph on $\mathbb CP^1$,
namely $f^{-1}(\mathbb RP^1)$. This
correspondence
is used for
example by S. Natanzon, B. Shapiro and A. Vainshtein to classify topologically generic real
rational maps
(see \cite{Shap1}
and \cite{Shap2}).
An other application of these graphs has been
exploited
by  A. Zvonkin
(\cite{Zvo}). He used these graphs to study the minimal degree of
$P^3-Q^2$,
where $P$ and $Q$ are complex polynomials of
degrees
$2k$ and $3k$,
respectively.
Following Zvonkin, Orevkov proposed in
 \cite{O3} a new way  to construct real algebraic
trigonal curves in rational geometrically ruled surfaces. Using
Orevkov's approach, we proved in \cite{Br1} the non-realizability of
some $\mathcal L$-schemes by real algebraic trigonal curves. All the
graphs considered in both these papers are special cases of real
rational graphs,  called \textit{real trigonal graphs} in \cite{Br1}. However, 
arguments used in
\cite{O3} are valid in the more general context of real rational graphs.

\begin{defi}
A root scheme is a k-uplet $((l_1,m_1),\ldots,(l_k,m_k))\in(\{p,q,r\}\times\mathbb N)^k$ with $k$ a natural number (here, $p,q$ and $r$ are symbols and do not stand for natural numbers).

\noindent A root scheme $((l_1,m_1),\ldots,(l_k,m_k))$ is realizable
by polynomials of degree $n$ if there exist two real polynomials in
one variable of degree $n$,
with no common roots, $P(X)$ and $Q(X)$ such that if
$x_1< x_2< \ldots < x_k$
are the real roots of $P,Q$ and $P+Q$, then
$l_i=p$ (resp., $q, r$) if $x_i$ is a root of $P$ (resp., $Q, P+Q$) and
$m_i$ is the multiplicity of $x_i$.

\noindent The polynomials $P$,  $Q$ and $P+Q$ are said to realize the root scheme $((l_1,m_1),\ldots,(l_k,m_k))$.
\end{defi}
In a root scheme, we will abbreviate
a sequence $S$
repeated $u$ times
by $S^u$.

\begin{figure}[h]
      \centering
 \begin{tabular}{cc}
\includegraphics[height=2cm, angle=0]{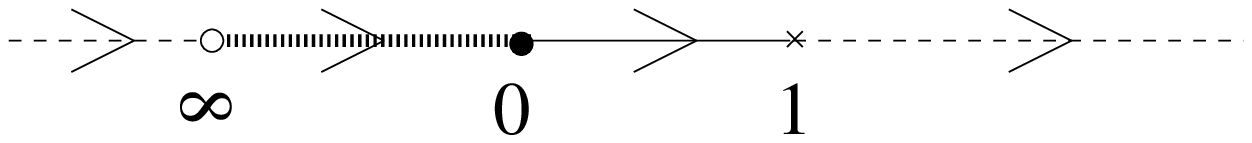}&
\includegraphics[height=2cm, angle=0]{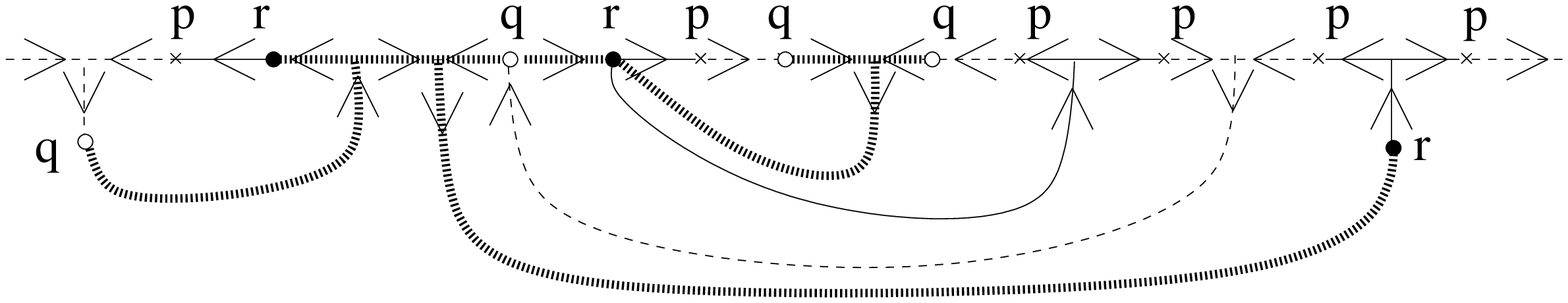}

\\ a)&b)
\end{tabular}
\caption{}
 \label{pregraph}
\end{figure}

From now on, let $RS$ be a root scheme and suppose that $RS$ is
realized by $P$, $Q$ and $P+Q$ of degree $n$.
Put $R(X)=P(X)+Q(X)$ and consider the rational function
$f(X)=\frac{R(X)}{Q(X)}=\frac{P(X)}{Q(X)}+1$. Color and orient
$\mathbb RP^1$ as depicted in Figure \ref{pregraph}a). Let $\Gamma$ be
 $f^{-1}(\mathbb RP^1)$ with the
coloring and the orientation induced by those chosen on $\mathbb
RP^1$.  Then, $\Gamma$ is a colored and oriented graph on $\mathbb
CP^1$, invariant under the action of the complex conjugation. 
The colored and oriented graph on $\mathbb RP^1$ obtained as the intersection of $\Gamma$ and $\mathbb RP^1$ can clearly be
extracted from $RS$.

\begin{defi}
The colored and oriented graph on $\mathbb RP^1$ constructed above is called the real graph associated to $RS$.
\end{defi}
\begin{defi}
Let $\Gamma$ be a graph
on $\mathbb CP^1$ invariant under the action
of the complex conjugation and $\pi:\Gamma\to\mathbb RP^1$ a
continuous map.
 Then the coloring and
orientation of $\mathbb RP^1$ shown in Figure \ref{pregraph}a) defines
a coloring and an orientation of $\Gamma$ via $\pi$.

\noindent The graph $\Gamma$ equipped with this coloring and this
orientation is called a real rational graph if 
\begin{itemize}
\item any vertex of $\Gamma$ has an even valence,
\item  any connected
component $D$ of $\mathbb CP^1\setminus\Gamma$ is homeomorphic to an open disk,
\item for any connected
component $D$ of $\mathbb CP^1\setminus\Gamma$, one has $\pi_{|\partial D}$ is a covering of $\mathbb RP^1$ of degree $d_D$.
\end{itemize}
The sum of the degrees $d_D$ for all connected
component $D$ of $\{Im(z)>0\}\setminus\mathbb RP^1$ of is called the degree of $\Gamma$.
\end{defi}
The importance of real rational graphs is given by the following
proposition.
\begin{prop}\label{constr poly}
Let $RS$ be a root scheme and $G$ its real graph. Then $RS$ is realizable by polynomials of degree $n$ if and only if there exists a real rational graph $\Gamma$ of degree $n$ such that $\Gamma\cap \mathbb RP^1=G$.
\end{prop}
\textit{Proof. }The proof used in \cite{O3} to construct real trigonal curves in $\Sigma_n$ can be used here in the same way.\findemo

For example, the root scheme $((p,1),(r,1),(q,2),(r,3),(p,1),(q,1),(q,1),(p,1),(p,1),(p,1),(p,1))$ is realizable by polynomials of degree $6$ as it is depicted on Figure \ref{pregraph}b).
\begin{lemma}\label{multiple roots}
Let $RS=((l_1,m_1),\ldots,(l_k,m_k))$ a root scheme such that there exist $i$ and $s$ such that $\forall j\in\{i,\ldots,i+s\}, l_j=l_i$.
Define the root scheme $RS'=((l_1',m_1'),\ldots,(l_{k-s}',m_{k-s}'))$ by
\begin{itemize}
\item $(l_t',m_t')=(l_t,m_t)$ for $t< i$,
\item $(l_i',m_i')=(l_i,m_i+\ldots+m_{i+s})$,
\item $(l_t',m_t')=(l_{t-s},m_{t-s})$ for $t> i+s$,
\end{itemize}
\noindent Then $RS$ is realizable by polynomials of degree $n$ if and only if $RS'$ is realizable by polynomials of degree $n$.
\end{lemma}
\textit{Proof. }Straightforward.\findemo

\section{Construction of reducible curves with a deep tangency point}\label{Constr tangency}
Let us define the root schemes $RS_n$ by
\begin{itemize}
\item $\left( (p,n),\left[ (r,1),(q,1)^2,(r,1) \right]^k,
(r,1),(q,1),(p,1),(r,1),
 \left[ (p,1),(r,1)^2,(p,1) \right]^k,(q,1)^n \right)$ if $n=2k+1$,
\item $\left( (p,n),\left[ (q,1),(r,1)^2,(q,1) \right]^{k},
(r,1),(q,1),(p,1),(r,1)^2,(p,1),
 \left[ (p,1),(r,1)^2,(p,1) \right]^{k},(q,1)^n \right)$ if $n=2k+2$.
\end{itemize}
\begin{prop}
For all $n$ in $\mathbb N^*$, the root scheme $RS_n$ is realizable by polynomials of degree $2n$.
\end{prop}
\textit{Proof. }According to lemma \ref{multiple roots}, one can substitute $(p,1)^n$ instead of $(p,n)$ in $RS_n$, and according to proposition \ref{constr poly}, one has just to construct a rational graph on $\mathbb CP^1$ with a real part corresponding to the real graph of $RS_n$. We will prove it by induction on $n$. All the pictures will represent the half $\{Im(z)\le 0\}$ of $\mathbb CP^1$.

\begin{figure}[h]
      \centering
 \begin{tabular}{cc}
\includegraphics[height=1.9cm, angle=0]{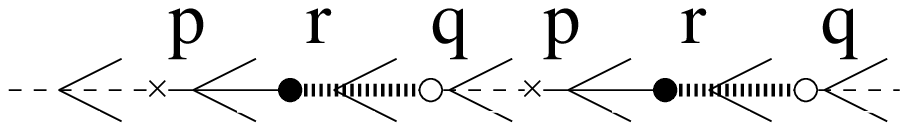}&
\includegraphics[height=1.9cm, angle=0]{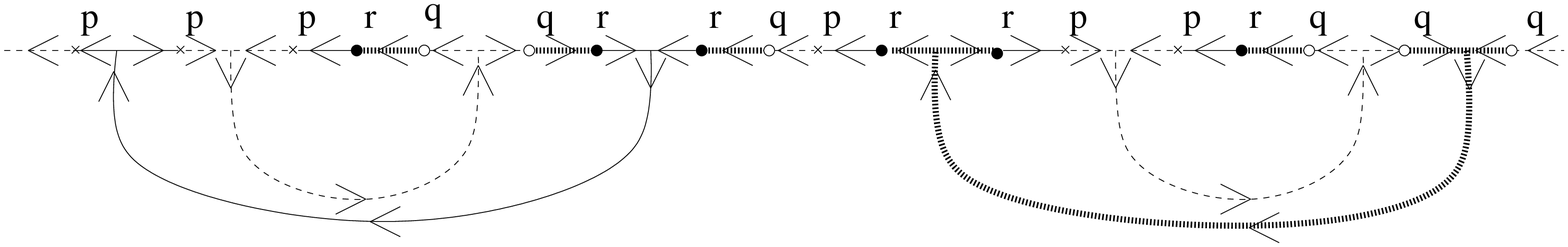}

\\ a)&b)
\end{tabular}
\caption{}
 \label{constr rational graph}
\end{figure}

The rational graph corresponding to $RS_1$ is depicted on Figure \ref{constr rational graph}a).
 \begin{figure}[h]
      \centering
 \begin{tabular}{cc}
\includegraphics[height=3.5cm, angle=0]{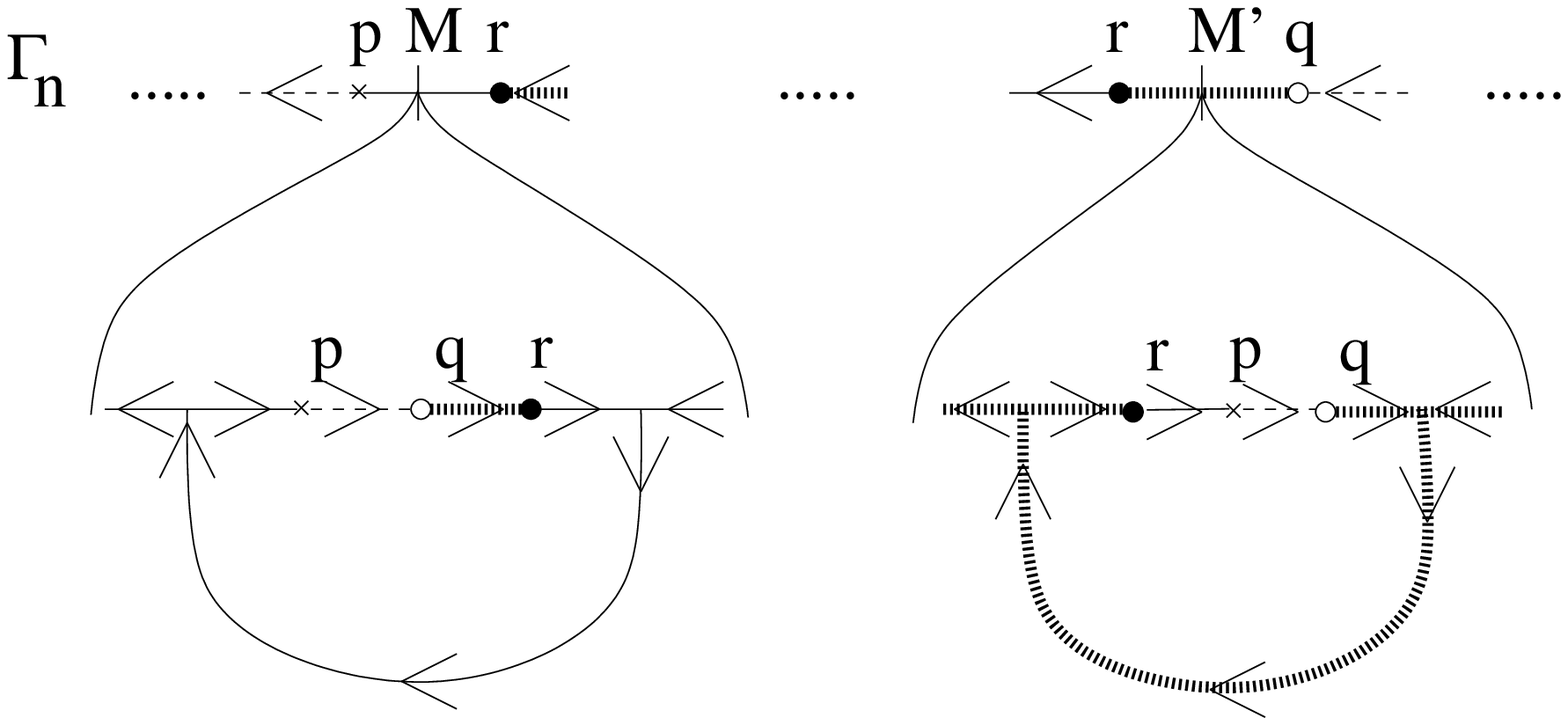}&
\includegraphics[height=3.5cm, angle=0]{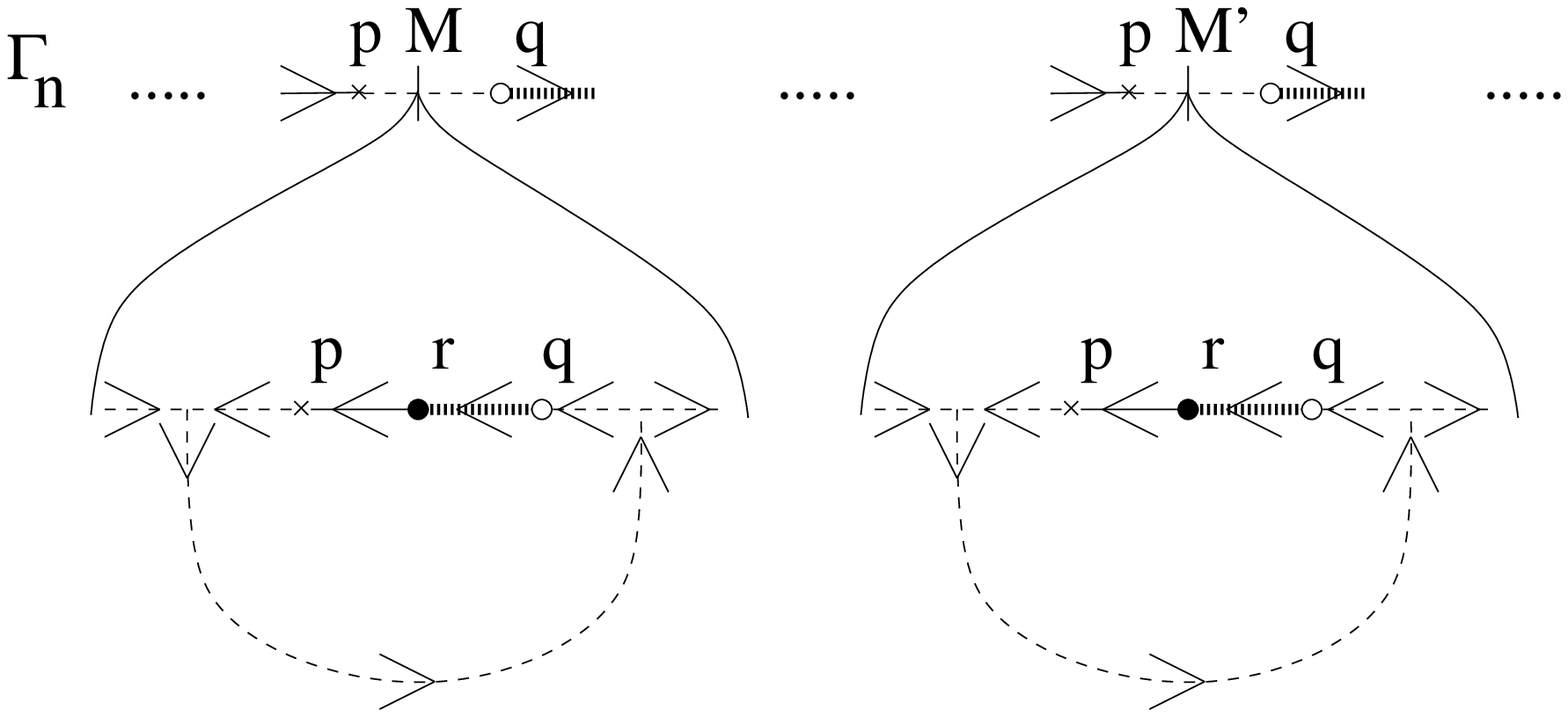}
\\ a)&b)
\end{tabular}
\caption{}
 \label{constr rational graph2}
\end{figure}

Suppose that a rational graph $\Gamma_n$ corresponding to $RS_n$ is constructed.

\noindent \textit{Consider, first, the case $n = 2k + 1$.}

Let $M$ be a point on $\Gamma_n\cap\mathbb RP^1$ between the
most right
point corresponding to $p$ in $RS_n$ and the
most left
point
corresponding to $r$ in $RS_n$. Then cut $\Gamma_n$ at $M$ and glue
the piece of graph depicted in Figure \ref{constr rational graph2}a).

 Let $M'$ be a point on $\Gamma_n\cap\mathbb RP^1$ between the more right point corresponding to $r$ in $RS_n$ and the more left point corresponding to $q$ in $RS_n$. Then cut $\Gamma_n$ at $M'$ and glue the piece of graph depicted in Figure \ref{constr rational graph2}a).

\noindent \textit{Consider now the case $n = 2k$. }

 Let $M$ be a point on $\Gamma_n\cap\mathbb RP^1$ between the more right point corresponding to $p$ in $RS_n$ and the more left point corresponding to $q$ in $RS_n$. Then cut $\Gamma_n$ at $M$ and glue the piece of graph depicted in Figure \ref{constr rational graph2}b).

 Let $M'$ be a point on $\Gamma_n\cap\mathbb RP^1$ between the more right point corresponding to $p$ in $RS_n$ and the more left point corresponding to $q$ in $RS_n$. Then cut $\Gamma_n$ at $M'$ and glue the piece of graph depicted in Figure \ref{constr rational graph2}b).

\vspace{2ex}
For example, $\Gamma_3$ is depicted in Figure \ref{constr rational graph}b). According to proposition \ref{constr poly}, the rational graphs $\Gamma_n$ ensure the realizability of the root schemes $RS_n$ by polynomials of degree $2n$.\findemo

\begin{cor}\label{poly}
For all $n$ in $\mathbb N^*$, there exists three real polynomials  $a_1(X)$, $a_2(X)$ and $b(X)$ of degree $n$ such that
\begin{itemize}
\item all the roots of $a_1$, $a_2$, $b$ and $a_1b+a2$ are real,
\item all the roots of $a_2$ and $a_1b+a2$ are smaller than the roots of $b$.
\end{itemize}
\end{cor}
\textit{Proof. }Let $P(X)$, $Q(X)$ and $R(X)=P(X)+Q(X)$ three polynomials of degree $2n$ realizing the root scheme $RS_n$. Then
\begin{itemize}
\item $Q(X)=\prod_{i=1}^{2n}(X-y_i)$ with $y_1<y_2<\ldots<y_{2n}$,
\item $P(X)=(X-\alpha)^n\prod_{i=1}^{n}(X-x_i)$ with $\alpha<x_1<x_2<\ldots<x_n<y_{n+1}$,
\item $R(X)=\prod_{i=1}^{2n}(X-z_i)$ with $z_1<z_2<\ldots<z_{2n}<y_{n+1}$.
\end{itemize}
Let us define $a_2(X)=(X+\alpha)^{2n}P\left(-\frac{1}{X+\alpha}\right)$, $A_1(X)=\prod_{i=1}^{n}(X-y_i)$, $a_1=(X+\alpha)^{n}A_1\left(-\frac{1}{X+\alpha}\right)$, $B(X)=\prod_{i=n+1}^{2n}(X-y_i)$ and $b=(X+\alpha)^{n}B\left(-\frac{1}{X+\alpha}\right)$.

\noindent As $a_1b+a2=(X+\alpha)^{2n}R\left(-\frac{1}{X+\alpha}\right)$, the corollary follows from the definition of $P$, $Q$ and $R$.\findemo

\vspace{4ex}
\noindent Now we are able to prove proposition \ref{C_n}.

\begin{figure}[h]
      \centering
 \begin{tabular}{cc}
\includegraphics[height=3cm, angle=0]{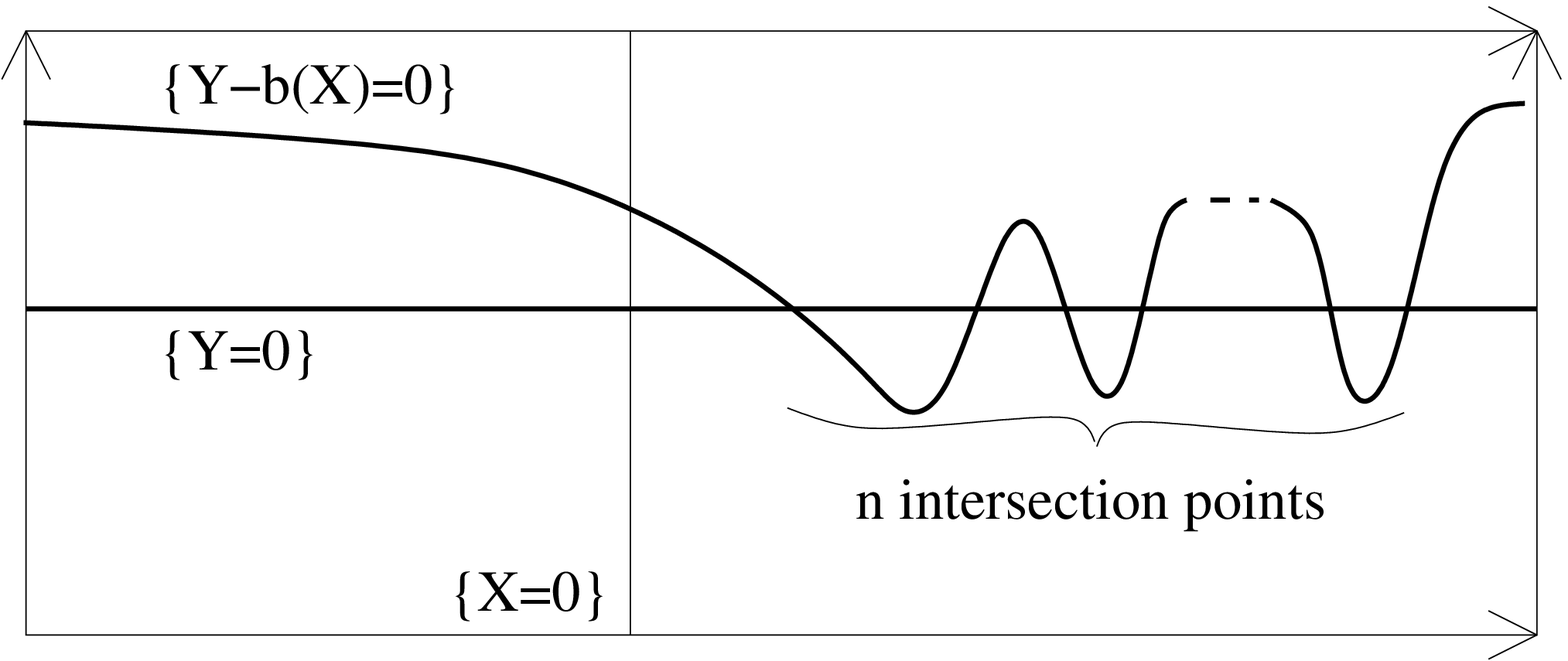}&
\includegraphics[height=3cm, angle=0]{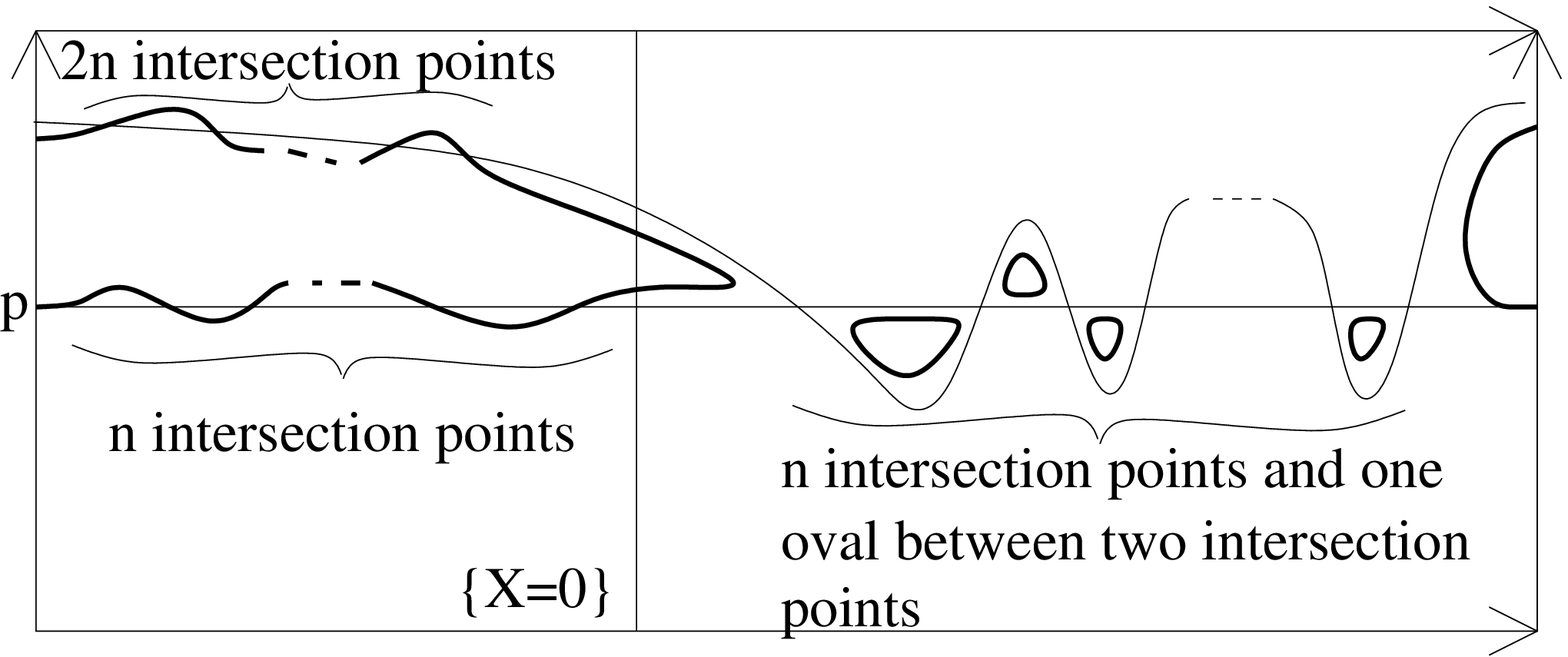}
\\ a)&b)
\end{tabular}
\caption{}
 \label{constr tang}
\end{figure}

\noindent \textit{Proof of proposition \ref{C_n}. }We construct
here explicitly only curves in $\Sigma_{2k}$. The construction of curves
in $\Sigma_{2k+1}$ is done in the same way.
Let us fix an even $n\ge 1$
and consider the polynomials $a_1(X)$, $a_2(X)$,
and $b(X)$ of degree $n$ constructed in corollary \ref{poly}.
Multiplying these three polynomials by $-1$ and performing a linear change
of
coordinates
if necessary, we can assume that the leading coefficient
of $b$ is positive,
all the roots of $b$ are
positive,
and
all the roots of $a_2$ and
$a_1b + a_2$
are
negative.
Then,
the
curve $Y(Y-b(X))$ in $\Sigma_n$ is depicted
in Figure \ref{constr tang}a).

\noindent For $t$ small enough and of suitable sign, the relative
positions of the curves $D_n(X,Y)=Y(Y-b(X))+t(a_1(X)Y+a_2(X))$,
$\{Y=0\}$ and $\{Y-b(X)=0\}$
are
as depicted in Figure  \ref{constr
  tang}b), where $p$ is a tangency point of order $n$ of $D_n$ and the
axis $\{Y=0\}$. Indeed, the definition of $a_1(X)$, $a_2(X)$,
and
$b(X)$  exactly means that the intersection points of $a_1(X)Y+a_2(X)$ and $Y(Y-b(X))$ have negative abscissa. Perturbing all the double points of $D_n(X,Y)(Y-b(X))$ in order to have the maximal number of ovals and keeping the tangency point of order $n$ with the axis $\{Y=0\}$, we obtain a curve $C_n$ with $n$ even, satisfying the conditions of proposition \ref{C_n}.\findemo

\section{Applications to real algebraic surfaces}\label{surfaces}
Here we recall and follow the notations proposed in \cite{B1}. We consider only $\mathbb Z/2\mathbb Z$-homology.

\noindent Let $d,i,k$ and $n$ be integers, $\beta_i(\mathbb RX^n_d)$
the $i^{th}$ Betti number of the real part of a nonsingular real hypersurface
of degree $d$ in $\mathbb CP^n$,
and $\beta_i(\mathbb RY^n_{2k})$ the $i^{th}$ Betti number
of the real part (for some real structure) of a double covering of $\mathbb CP^n$ branched
along a nonsingular real hypersurface of degree $2k$.

\noindent An interesting
question
concerns
the asymptotic behavior of $max$
$\beta_i(\mathbb RX^n_d)$ and $max$ $\beta_i(\mathbb RY^n_{2k})$ when
$d$ and $k$ go to infinity.
In \cite{B1}, Bihan has showed that there exist two sequences
of real numbers $(\zeta_{i,n})_{i,n\in\mathbb N^2}$
and $(\delta_{i,n})_{i,n\in\mathbb N^2}$
such that
\begin{center}
\begin{tabular}{ccccccccc}
$max\textrm{ }\beta_i(\mathbb RX^n_d)$&$\sim$&$\zeta_{i,n}d^n$&&
and  &&$max\textrm{ }\beta_i(\mathbb RY^n_{2k})$&$\sim$&$\delta_{i,n}k^n.$
\\&$\scriptstyle{d\to\infty}$&&&&&&$\scriptstyle{k\to\infty}$
\end{tabular}
\end{center}
The exact value of the numbers $\zeta_{i,n}$ and $\delta_{i,n}$ are
known only for small $n$. The following equalities are well known (see \cite{B1}, \cite{I1}).
$$\delta_{0,0}=2,\textrm{ }\zeta_{0,1}=\delta_{0,1}=\delta_{1,1}=1, \textrm{ }\zeta_{0,2}=\zeta_{1,2}=\frac{1}{2},$$
$$\delta_{0,2}\le\frac{7}{4},\textrm{ }\delta_{1,2}\le\frac{7}{2}, \textrm{ }\zeta_{0,3}\le\frac{5}{12}\textrm{ and }\zeta_{1,3}\le\frac{5}{6}.$$

The upper bounds are classical and are obtained using the Harnack  and
Comessati-Petrovsky-Oleinik inequalities. Lower bounds for
$\delta_{0,2}$ and $\delta_{1,2}$ are directly related to the
asymptotically maximal number of even ovals of a curve of even degree
in $\mathbb RP^2$, and before the results of the present paper, the
best known lower bounds for these two numbers were, respectively,
$\frac{27}{16}$ and $\frac{27}{8}$ (see \cite{I1}). In \cite{B1},
Bihan has constructed nonsingular real algebraic surfaces in $\mathbb RP^3$ with Betti numbers related to $\delta_{0,2}$ and $\delta_{1,2}$.
\begin{thm}[Bihan]
One has $\frac{\delta_{0,2}}{6}+\frac{1}{12}\le\zeta_{0,3}
\textrm{ and }\frac{\delta_{1,2}}{6}+\frac{1}{6}\le\zeta_{1,3}. $
\end{thm}

\noindent Theorem \ref{max even} gives as immediate corollaries the exact values of
$\delta_{0,2}$ and $\delta_{1,2}$ and improves the known lower bounds
for $\zeta_{0,3}$ and $\zeta_{1,3}$.
\begin{prop}
One has $\delta_{0,2}=\frac{7}{4}$ and $\delta_{1,2}=\frac{7}{2}.$
\end{prop}
\begin{cor}
One has $\frac{9}{24}\le\zeta_{0,3}\le\frac{5}{12}
\textrm{ and }\frac{9}{12}\le\zeta_{1,3}\le\frac{5}{6}. $
\end{cor}

\small
\def\rightmark{\em Bibliography}
\addcontentsline{toc}{section}{References}
\bibliographystyle{alpha}
\bibliography{/home/erwan/These/Biblio.bib}

\hspace*{2 ex}
\textbf{\\Erwan Brugall\'e}\\
Laboratoire Emile Picard\\
Université Paul Sabatier\\
UFR MIG\\
118 route de Narbonne\\
31 000 Toulouse\\
FRANCE\\
\\
E-mail : brugalle@picard.ups-tlse.fr
\end{document}